\documentclass[11pt,oneside,english]{amsart}
\usepackage{geometry}
\usepackage{amssymb}

\makeatletter
 \theoremstyle{plain}    
 \newtheorem{thm}{Theorem}[section]
 \numberwithin{equation}{section} 
 \numberwithin{figure}{section} 
 \theoremstyle{plain}
 \theoremstyle{remark}
 \newtheorem{rem}[thm]{Remark}
 \theoremstyle{definition}
 \newtheorem{defn}[thm]{Definition}
 \theoremstyle{plain}    
 \newtheorem{prop}[thm]{Proposition} 
 \theoremstyle{remark}    
 \newtheorem{case}{Case} 

\newcommand{\R}{\mathbb{R}}
\newcommand{\C}{\mathbb{C}}

\usepackage{babel}
\makeatother
\begin{document}

\title{Volume-minimizing foliations on spheres}

\author{Fabiano Brito and David L. Johnson}

\address{Dpto. de Matem\'{a}tica, Instituto de Matem\'{a}tica e Estat\'{i}stica,
Universidade de S\={a}o Paulo, R. do Mat\={a}o 1010, S\={a}o Paulo-SP
05508-900, Brazil.}

\email{fabiano@ime.usp.br}

\address{Department of Mathematics, Lehigh University, 14 E. Packer Ave.,
Bethlehem, PA  18015, USA}

\email{david.johnson@lehigh.edu}

\subjclass{53C12, 53C38}

\thanks{The second author was supported during this research by grants from
the Universidade de S\={a}o Paulo, FAPESP Proc.1999/02684-5, and
Lehigh University, and thanks those institutions for enabling the
collaboration involved in this work.}

\begin{abstract}
The volume of a $k$-dimensional foliation $\mathcal{F}$ in a Riemannian
manifold $M^{n}$ is defined as the mass of image of the Gauss map,
which is a map from $M$ to the Grassmann bundle of $k$-planes in
the tangent bundle. Generalizing the construction by Gluck and Ziller
in \cite{GZ}, {}``singular'' foliations by $3$-spheres are constructed
on round spheres $S^{4n+3}$, as well as a singular foliation by $7$-spheres
on $S^{15}$, which minimize volume within their respective relative
homology classes. These singular examples provide lower bounds for
volumes of regular $3$-dimensional foliations of $S^{4n+3}$ and
regular $7$-dimensional foliations of $S^{15}$.
\end{abstract}
\maketitle

\section{Introduction}

In \cite{GZ}, Herman Gluck and Wolfgang Ziller asked which foliations
were {}``best-organized'', in that an energy functional they called
the \emph{volume} was minimized. The \emph{volume} of a foliation
is the mass of the image of the Gauss map, which in the case of a
one-dimensional foliation is the mass of the unit tangent flow field
in $T_{1}(M)$. 

They were able to show that the standard one-dimensional foliation
(or \emph{flow}, in their terminology) of $S^{3}$ by the fibers of
the Hopf fibration $S^{3}\rightarrow S^{2}$ minimized volume among
all foliations of the round $S^{3}$. Their method of proof, involving
calibrations, did not generalize, however. 

It is not the case that even the most obvious generalization of Gluck
and Ziller's example to higher dimensions, the Hopf fibration $S^{5}\rightarrow\C\mathbb{P}^{2}$,
is volume-minimizing \cite{Jo}. Sharon Pedersen showed in her thesis
that there was a foliation of $S^{5}$ with much less volume than
the Hopf fibration, although her example is singular \cite{P}. It
may well be that the volume-minimizing one-dimensional foliations
on $S^{5}$ is be singular, although it is not clear whether Pedersen's
example is that minimizer. Gluck and Ziller did describe a {}``singular
foliation'' on $S^{2n+1}$ that minimizes the volume functional,
but their singular minimum is of a different sort than Pedersen's.
Pedersen's foliation is a smooth foliation on all but one point in
$S^{5}$, and is a limit of smooth foliations, while Gluck and Ziller's
example is not homologous to a foliation except on $S^{3}$.

There is, then, something peculiar about the Hopf fibration on $S^{3}$
which enables the calibration argument that Gluck and Ziller used
to show the minimization of the volume of that foliation, beyond the
evident geometric properties for the Hopf fibrations in general. 

In this article we expand the method used by Gluck and Ziller to 3-dimensional
foliations of $S^{4n+3}$ and 7-dimensional foliations of $S^{15}$.
What we find is that the generic situation Gluck and Ziller described
for flows on $S^{2n+1}$ holds; that is, there are singular foliations
which minimize volume in these cases, but that it does not appear
that the Hopf fibrations will minimize volume.

\section{Definitions and the minimization question}

The original question considered by Gluck and Ziller in \cite{GZ},
extended by a number of authors, is to find the dimension-$k$ foliation
$\mathcal{F}$ on a compact Riemannian manifold $M$, considered as
a section \[
\sigma_{\mathcal{F}}:M\rightarrow G_{o}(k,M)\]
of the bundle of oriented $k$-planes tangent to $M$, which is {}``most
efficient'' or {}``best-organized'' in that its \emph{volume} is
minimized, where the \emph{volume} is defined as the Hausdorff $n$-dimensional
measure of the image $\sigma_{\mathcal{F}}(M)\subset G_{o}(k,M)$,
where the Grassmann bundle has a natural Sasaki metric induced from
the original metric on $M$. Volume-minimization should be considered
within each homology class of foliations, and it is possible for one
homology class to admit a smooth minimizer, but for others to have
no smooth minimizer. 

\begin{rem}
It may seem more appropriate to consider homotopy classes of such
foliations rather than homology classes, but a simple construction
shows that two homotopy classes of one-dimensional foliations on $S^{3}$
can be constructed (within one homology class, of course), one of
which has a smooth volume-minimizer, but the other does not, since
there is a sequence within the one homotopy class whose volume converges
to the minimum of the other class. Since the only foliations achieving
that minimum are within the first homotopy class (see, for example,
\cite{GZ}), there can be no smooth minimizer within the first. 
\end{rem}
As mentioned in the Introduction, Gluck and Ziller showed that the
natural candidate, the fibers of the Hopf fibration from $S^{3}$
to $S^{2}$, is volume-minimizing among all (smooth) one-dimensional
foliations on the (round) 3-sphere. 

Several authors \cite{Jo,P} showed that this natural candidate volume-minimizer
did not extend even to the next simplest case of the fibers of the
Hopf fibration $S^{5}\rightarrow\C\mathbb{P}^{2}$. Pedersen's example,
in particular, is singular in the sense that there is one point of
$S^{5}$ which must be removed in order for her example to be a smooth
foliation. It is the case, however, that Pedersen's example is the
limit of smooth foliations (it is the limit of the sequence of geodesic
flows stretching away from one pole towards the other, applied to
any smooth one-dimensional foliation). 

Because of Pedersen's example, it seems necessary to consider singular
foliations in general.

\begin{defn}
An \emph{oriented} \emph{singular $k$-dimensional distribution} on
a manifold $M$ is defined as an $n$-dimensional rectifiable current
$D\subset G_{o}(k,M)$ of the bundle of $k$-dimensional subspaces
of $T_{*}(M)$, so that on an open dense subset $U\subset M$, $\left.D\right|_{\pi^{-1}(U)}$
is a smooth, $k$-dimensional distribution on $U$, that is, a smooth
cross-section of $G(k,U)\rightarrow U$ (resp., $G_{o}(k,U)\rightarrow U$).
The distribution $D$ is \emph{integrable}, or is a \emph{singular
foliation}, if $\left.D\right|_{\pi^{-1}(U)}$ is integrable.

As an example, any unit vector field on a manifold $M$ with finitely
many singularities, each with finite index, is an oriented singular
foliation in this sense. Note that these currents need not be cycles,
in general; for example in the case of a unit vector field with some
point singularity of odd degree.

This notion of a singular foliation is similar to, but more general
than, that studied by the second-named author and Smith in \cite{JS}.
In that article, the singular sections of arbitrary vector bundles
that are considered are those in the weak closure of the space of
smooth sections. Many of the singular foliations considered here are
not in the closure of the space of smooth sections, by topological
considerations.
\end{defn}

\section{1-dimensional singular foliations of $S^{2n+1}$}

\subsection{The calibration}

The bundle of oriented 1-planes tangent to $S^{2n+1}$, the unit tangent
bundle $T_{1}(S^{2n+1})$, is isomorphic to the flag manifold of oriented
lines in oriented 2-planes in $\R^{2n+2}$, which is the Stieffel
manifold of $2$-frames in $\R^{2n+2}$. 

This gives rise to the following diagram:\[
\begin{array}{ccc}
T_{1}(S^{2n+1}) & \longrightarrow & F_{o}(1,2,\R^{2n+2})\\
h\Big{\uparrow}\pi\Big{\downarrow} &  & \pi\Big{\downarrow}\\
S^{7} &  & G_{o}(2,\R^{2n+2})\end{array}.\]
$G_{o}(2,2n+2)$ has two universal bundles, the universal 2-plane
bundle $U(2,2n+2)$ and the dual $2n$-plane bundle $V(2n,2n+2)$,
defined by \begin{eqnarray*}
U(2,2n+2) & := & \cup_{x\in G_{o}(2,2n+2)}x\\
V(2n,2n+2) & := & \cup_{x\in G_{o}(2,2n+2)}x^{\perp}.\end{eqnarray*}
The respective Euler classes $E(U)$ and $E(V)$ satisfy $E(U)\cup E(V)=0$
in $H^{2n}(G_{o}(2,2n+2))$, since $U\oplus V$ is trivial. In particular,
if $\omega$ is the universal connection on $U(2,2n+2)$ defined by
Narasimhan and Ramanan (cf. \cite{NR}), and $\omega^{*}$ is the
{}``dual'' connection on $V(2n,2n+2)$, then the associated Euler
forms, $e(\Omega)$ and $e(\Omega^{*})$, satisfy $e(\Omega)\wedge e(\Omega^{*})=0$.
Consider the form \[
\Phi:=C\, Te(\omega)\wedge e(\Omega^{*}),\]
which is well-defined on $F_{o}(1,2,\R^{2n+2})$ since that is the
frame bundle $FU(2,2n+2)$ of oriented orthonormal frames on $U(2,2n+2)$,
which is an $SO(2)$-principal bundle. Here, $Te(\omega)$ is the
transgressive Chern-Simons form corresponding to the Euler form $e(\Omega)$
of $U(2,2n+2)$ \cite{CS}. Because $d(Te(\omega))=e(\Omega)$, we
have that $d\Phi=0$. The constant $C$ is simply chosen so that the
comass of $\Phi$ is one. This is the same calibration defined in
\cite{GZ}.

\subsection{Calculations}

We will consider $G_{o}(2,2n+2)$ as $SO(2n+2)/SO(2)\times SO(2n)$,
and the principal bundle $FU(2,2n+2)$ as $SO(2n+2)/I_{2}\times SO(2n)$.
The universal connection $\omega$ on $FU(2,2n+2)$ can be defined
as the truncation of the restriction of the Maurer-Cartan form on
$o(2n+2)$, denoted $\mu=\left[\mu_{ij}\right]$, to the tangents
to $FU(2,2n+2)$. That is, the components of the connection $\omega_{ij}$
are defined for $i,j\in\{1,2\}$ by $\omega_{ij}(A)=A_{ij}$, for
any \[
A\in T_{*}(FU(2,2n+2),(U_{0},\left\{ e_{1},e_{2}\right\} ))=\left\{ \left.A\in o(2n+2)\right|A=\left[\begin{array}{cc}
R & S\\
-S^{t} & 0\end{array}\right],\, R\in o(2)\right\} ,\]
if $U_{0}=\R^{2}\times0\subset\R^{2n+2}$, with basis $\left\{ e_{1},e_{2}\right\} $.
By homogeneity, all calculations in $FU(2,2n+2)$ can be taken to
be at this point. 

The curvature $\Omega$ of this connection is given by $\Omega_{ij}(X,Y)=-\omega_{ij}([X,Y])$
for left-invariant vector fields that are horizontal at $U_{0}$,
that is, of the form $\left[\begin{array}{cc}
0 & S\\
-S^{t} & 0\end{array}\right]$. In terms of the Maurer-Cartan form, $\Omega_{ij}=+\sum_{k=2}^{2n+2}\mu_{ik}\wedge\mu_{jk}$,
for $i,j\in\{1,2\}$. 

Similarly, the connection $\omega^{*}$ on the dual principal bundle
$FV(2n,2n+2)=SO(2n+2)/SO(2)\times I$ at $U_{0}=\R^{2}\times0\subset\R^{2n+2}$
is the restriction of the same Maurer-Cartan form $\mu$ to the other
block, and the curvature $\Omega_{kl}^{*}=\sum_{i=1}^{2}\mu_{ik}\wedge\mu_{il}$,
for $k,l\in\left\{ 3,\ldots,2n+2\right\} $. Either of the tangent
spaces to these principal bundles can be canonically embedded into
the tangent space $o(2n+2)$ of $SO(2n+2)$ at the identity.

The Euler form $e(\Omega)$ of $FU(2,2n+2)$ is the form \begin{eqnarray*}
e(\Omega) & := & \frac{1}{2\pi}\left(\Omega_{12}\right)\\
 & = & \frac{1}{2\pi^{2}}\left(\mu_{1k}\wedge\mu_{2k}\right),\end{eqnarray*}
where the sum is taken over all $k\in\left\{ 3,\ldots,2n+2\right\} $.
Dually, the Euler form $e(\Omega^{*})$ of $FV(2n,2n+2)$ is the form
\begin{eqnarray*}
e(\Omega^{*}) & := & C\left(\sum_{\sigma\in S_{2n}}(-1)^{\sigma}\Omega_{\sigma(3)\sigma(4)}\wedge\cdots\wedge\Omega_{\sigma(2n+1)\sigma(2n+2)}\right)\\
 & = & C\left(\sum_{\sigma\in S_{2n},i_{1},\ldots,i_{n}}(-1)^{\sigma}\mu_{\sigma(3)i_{1}}\wedge\mu_{\sigma(4)i_{1}}\wedge\cdots\wedge\mu_{\sigma(2n+1)i_{n}}\wedge\mu_{\sigma(2n+2)i_{n}}\right),\end{eqnarray*}
where the sum is taken over all $\sigma\in S_{2n}$as permutations
of $\{3,\ldots,2n+2\}$, $i_{1},\ldots,i_{n}\in\left\{ 1,2\right\} $,
and the constant depends just on the dimension.

\begin{prop}
$e(\Omega)\wedge e(\Omega^{*})\equiv0$.
\end{prop}
\begin{proof}
Each monomial in this product is of the form \[
\mu_{1k}\wedge\mu_{2k}\wedge\mu_{3i_{1}}\wedge\mu_{4i_{1}}\wedge\cdots\wedge\mu_{(2n+1)i_{n}}\wedge\mu_{(2n+2)i_{n}}\]
 or a permutation thereof. $k$ can be in $3,\ldots,2n+2$. No matter
what $k$ is, since $i_{1},\ldots i_{n}$ are either 1 or 2, then
this form must be 0. 
\end{proof}
Thus, the form $\Phi:=C\, Te(\omega)\wedge e(\Omega^{*})$ is indeed
closed. 

It remains to find the maximum of $\Phi(W)$ for $2n+1$-planes $W$
in the total space of $\pi\: F(1,2,\R^{2n+2})\rightarrow G(2,2n+2)$.

Certainly the vertical direction will be a maximum for $Te(\omega),$
which is (up to scale) exactly the volume form of the fibers. Thus
the maximum is achieved only when one direction of the $(2n+1)$-plane
is vertical. 

It is interesting to note that, since the maximum of $\Phi$ must
necessarily have a vertical direction at each point, any current calibrated
by $\Phi$ must be a contained in a union of fibers of the projection
$\pi:F(1,2,\R^{2n+2})\rightarrow G(2,2n+2)$, so must be of the form
$\pi^{-1}(M)\cap U$ for some current $M\subset G(2,2n+2)$. Since,
for $W\in G(2,2n+2)$, the preimage \[
\pi^{-1}(W)=\left\{ x\left|x\in W,\,|x|=1\right.\right\} =\left\{ \{ e_{1},e_{2}\}\left|\{ e_{1},e_{2}\}\,\textrm{is a basis of}\, W\right.\right\} \]
 is, as a subset of $T_{1}(S^{2n+1})$, the unit velocity field of
the great circle $S^{2n+1}\cap W$ with orientation determined by
$W$. In terms of the foliations determined by these calibrated currents,
they must then consist of arcs of great circles, and must be great
circle foliations if they are regular. 

To see what currents $\Phi$ calibrates, we now need only find those
$2n$-plane directions maximizing $e(\Omega^{*})$.

Since \begin{eqnarray*}
e(\Omega^{*}) & := & C\left(\sum_{\sigma\in S_{2n}}(-1)^{\sigma}\Omega_{\sigma(3)\sigma(4)}\wedge\cdots\wedge\Omega_{\sigma(2n+1)\sigma(2n+2)}\right)\\
 & = & C\left(\sum_{\sigma\in S_{2n},i_{1},\ldots i_{n}}(-1)^{\sigma}\mu_{\sigma(3)i_{1}}\wedge\mu_{\sigma(4)i_{1}}\wedge\cdots\wedge\mu_{\sigma(2n+1)i_{n}}\wedge\mu_{\sigma(2n+2)i_{n}}\right),\end{eqnarray*}
if $E_{ij}$ is the basis of tangent vectors dual to $\mu_{ij}$,
for any fixed permutation $\sigma\in S_{2n}$, \[
e(\Omega^{*})(E_{1\sigma(3)},E_{1\sigma(4)},\ldots,E_{1\sigma(2n+2)})=(-1)^{\sigma}(2n)!C=e(\Omega^{*})(E_{2\sigma(3)},\ldots,E_{2\sigma(2n+2)}).\]
 It is straightforward to see that, if $i_{j_{1}}\neq i_{j_{2}}$,
then some permutations in the sum will evaluate to 0, so that\[
\left|e(\Omega^{*})(E_{i_{1}\sigma(3)},E_{i_{2}\sigma(4)},\ldots,E_{i_{2n}\sigma(2n+2)})\right|<(2n)!C.\]
 Finally, if $\{ i_{1},\ldots,i_{2n}\}$ does not have at least $n$
pairs of values, or if $\{ k_{1},\ldots,k_{2n}\}$ does not consist
of some permutation of $\{3,\ldots,2n+2\},$ then $e(\Omega^{*})(E_{i_{1}k_{1}},\ldots,E_{i_{2n}k_{2n}})=0$. 

For any decomposable, unit $\xi\in\Lambda_{2n}(G(2,2n+2),W_{0})$
which is tangent to the variety $G(2,2n+2)$ at $W_{0}$, \[
\xi=\sum_{i_{1},\ldots,i_{2n},k_{1}\leq\cdots\leq k_{2n}}\xi_{i_{1},\ldots,i_{2n},k_{1},\ldots,k_{2n}}E_{i_{1}k_{1}}\wedge\cdots\wedge E_{i_{2n}k_{2n}}.\]
Since $\xi$ is decomposable, $\xi$ satisfies the Pl\"{u}cker condition
$\xi\wedge\xi=0$, implying that, in particular (restricting to the
case where $\{ k_{1},\ldots,k_{2n}\}=\{3,\ldots,2n+2\}$ since otherwise
$e(\Omega^{*})=0$), and denoting $\xi_{i_{1},\ldots,i_{2n},3,\ldots,(2n+2)}$
by $\xi_{i_{1},\ldots,i_{2n}}$, \[
\xi_{1,\ldots,1}\xi_{2,\ldots,2}-\xi_{2,1,\ldots,1}\xi_{1,2,\ldots,2}-\xi_{1,2,1,\ldots,1}\xi_{2,1,2,\ldots,2}+\cdots=0,\]
and similarly for all other such combinations. Thus, \begin{eqnarray*}
(\xi_{1,\ldots,1}+\xi_{2,\ldots,2})^{2} & = & \xi_{1,\ldots,1}^{2}+\xi_{2,\ldots,2}^{2}+2\xi_{1,\ldots,1}\xi_{2,\ldots,2}\\
 & = & \xi_{1,\ldots,1}^{2}+\xi_{2,\ldots,2}^{2}+2\xi_{2,1,\ldots,1}\xi_{1,2,\ldots,2}+2\xi_{1,2,1,\ldots,1}\xi_{2,1,2,\ldots,2}\pm\cdots\\
 & \leq & \xi_{1,1,1,1}^{2}+\xi_{2,2,2,2}^{2}+\xi_{2,1,\ldots,1}^{2}+\xi_{1,2,\ldots,2}^{2}+\xi_{1,2,1,\ldots,1}^{2}+\xi_{2,1,2,\ldots,2}^{2}+\cdots\\
 & \leq & 1,\end{eqnarray*}
since $\xi$ is a unit. Thus, on any such $\xi$, \[
e(\Omega^{*})(\xi)\leq(2n)!C,\]
the maximum being achieved on those $\xi$ so that $(\xi_{1,\ldots,1,3,\ldots,(2n+2)}+\xi_{2,\ldots,2,3,\ldots,(2n+2)})=1$
which have the proper orientation. Those $2n$-planes are, except
where $n=1$, \emph{not} those which are complex $2n$-planes in $T_{*}(G(2n,2n+2),W_{0})$
under some complex structure on that space induced from one of $\R^{2n+2}$
for which $W_{0}$ is complex. 

\begin{thm}
The standard foliation $H$ of $S^{3}$ by the fibers of the Hopf
fibration $S^{3}\rightarrow S^{2}$ for some complex structure on
$\R^{4}\supset S^{3}$ minimizes the volume of one-dimensional foliations
of $S^{3}$. The singular foliation $NS$ of $S^{2n+1}$, $n>1$ consisting
of all great circles through a pair of antipodal points with indices
$\pm1$ minimizes volume of all singular foliations on $S^{2n+1}$
with those singular points and indices, and provides a lower bound
for the volume of all one-dimensional oriented foliations of $S^{2n+1}$.
\end{thm}
\begin{rem}
The minimization of the Hopf fibration in the case $n=1$ is due to
Gluck and Ziller in \cite{GZ}. They also showed a bound on the minimum-volume
flow in higher dimensions by constructing a specific cycle in twice
the homology of a flow. The first-named author, along with P. Chac\'{o}n
and A. M. Naveira, in \cite{BCN}, showed that this bound is attained
by the specific singular foliation $NS$, and is a strict lower bound
for volumes of smooth foliations.  The notation $NS$ ({}``north-south'')
refers to the fact that this foliation is by longitude lines from
one pole to the other.
\end{rem}
\begin{proof}
~
\begin{case}
$n=1$
\end{case}
In the case $n=1$ any complex $2$-plane will maximize $e(\Omega^{*})$,
since $T_{*}(G(2,4),W_{0})$ is $\C^{2}$ and a real 2-plane $\xi$
in $\C^{2}$ is complex (for a given complex structure) if and only
if $<\xi,\alpha>+<\xi,\beta>=1$ for any orthogonal pair $\alpha,\beta$
of complex lines, where the inner product is the standard induced
inner product on $\Lambda_{2}(\C^{2})$ induced from the inner product
on $\C^{2}$ itself. Using coordinate planes and the standard complex
structure on $G(2,4)$ (which is as the projective variety in $\C\mathbb{P}^{3}$
defined by $z_{0}^{2}+z_{1}^{2}+z_{2}^{2}+z_{3}^{2}=0$), this condition
is equivalent to $(\xi_{1,1,3,4}+\xi_{2,2,3,4})=1.$ So any complex
submanifold $M\subset G(2,4)$ will be calibrated by $e(\Omega^{*})$. 

Not every such complex submanifold corresponds to a foliation of $S^{3}$,
however, not even a singular one. For any $W\in M$, the preimage
$\pi^{-1}(W)\subset F(1,2,4)=T_{1}(S^{3})$ corresponds to the image
in $T_{1}(S^{3})$ of the intersection of $S^{3}$ with the 2-plane
$W$ via the tangent map, a great circle on $S^{3}$. So, if $M$
(complex or not) corresponds to a smooth or singular foliation of
$S^{3}$, it is a foliation by great circles. If all such $W$ are
complex lines in $\R^{4}=\C^{2}$ for some complex structure on the
$\R^{4}$ in which $S^{3}$ is embedded, then all of these great circles
are disjoint, $M$ is the standard embedding of $\C\mathbb{P}^{1}$
in $G_{0}(2,4)$, and the foliation is a Hopf fibration, and the corresponding
curve $M$ in $G_{0}(2,4)\subset\C\mathbb{P}^{3}$ is defined by $z_{0}=iz_{1}$
in addition to $z_{0}^{2}+z_{1}^{2}+z_{2}^{2}+z_{3}^{2}=0$. Other
complex submanifolds of $G_{0}(2,4)$ do not correspond to even a
singular foliation of $S^{3}$. For example, the curve $z_{0}=0$,
which is also a hyperplane section of $G_{0}(2,4)$ and which is $G_{0}(2,3)\cong\C\mathbb{P}^{1}$,
lifts to $F_{0}(1,2,3)=T_{1}(S^{2})\subset T_{1}(S^{3})$, so does
not correspond to a section over a dense subset of $S^{3}$. 

The manifold $M=\left\{ \left.W\in G_{0}(2,4)\right|e_{1}\in W\right\} $,
which is dual to the previous submanifold, will also be calibrated
by $\Phi$, since at $W_{0}\in M$, with basis chosen so that $W_{0}=e_{1}\wedge e_{2}\in\Lambda_{2}(\R^{4})$,
the tangent plane satisfies $\xi_{2,2,3,4}=1$. The current $NS$
corresponding to a singular foliation will not be all of $\pi^{-1}(M)$,
since that will be a double of the singular foliation by all great
circles passing through $\pm e_{1}$. Instead, the current $NS$ is
formed from semicircular fibers of this bundle, from the fiber of
$T_{1}(S^{3})$ over $-e_{1}$ to that over $+e_{1}$. This current
would minimize volume over all singular foliations of $S^{3}$ with
two point singularities at $\pm e_{1}$, $-e_{1}$ having index $-1$
and $e_{1}$ having index $1$. The minimum volume of such singular
foliations is the same as that of the Hopf fibrations.

\begin{case}
$n>1$
\end{case}
For $n>1$ if $M$ is the manifold \[
M:=\left\{ \left.W\in G_{0}(2,2n+2)\right|e_{1}\in W\right\} :=\left\{ \left.e_{1}\wedge x\right|x\perp\{ e_{1}\},\,\left\Vert x\right\Vert =1\right\} ,\]
then $M\cong S^{2n}$ is not the space of complex $2n$-planes in
$\R^{2n+2}$ for any complex structure, and the corresponding {}``foliation''
on $S^{2n+1}$ will be singular. The tangent planes to $M$ at each
point clearly maximize the value of $e(\Omega^{*})$. Note also that,
in this case complex submanifolds of $G(2,2n+2)$ are \emph{not} calibrated
by $e(\Omega^{*})$. 

In general, this singular distribution will indeed be calibrated by
this form, so minimizes volume among all singular foliations with
the same singular set; in this case, an antipodal pair of singular
points, with indices $\pm1$ that are in each leaf of the singular
foliation. Since the current in $T_{1}(S^{2n+1})$ actually defined
by $M$ consists of the unit tangent field to oriented semi-circles,
longitudes, from $-e_{1}$ to $+e_{1}$ in $S^{2n+1}$, which has
as a 2-fold cover the submanifold $S^{2n}\times S^{1}=\pi^{-1}(M)\subset F_{0}(1,2,\R^{2n+2})=T_{1}(S^{2n+1})$,
the mass-minimization property of the calibration compares the mass
of this current, $NS$, to all other currents $S$ with the same boundary
(the two tangent fibers over $\pm e_{1}$, suitably oriented), which
are homologous in that $NS-S$ is a boundary. This can be easily extended
to all other currents with the same singular points and the same indices
at those singular points, since any such current can be modified within
the singular fibers to match the boundary of $NS$.

\end{proof}

\section{3-dimensional foliations of $S^{7}$}

\subsection{The calibration}

Note that the Grassmann bundle $G(3,S^{7})$ of oriented 3-planes
tangent to $S^{7}$ is isomorphic to the flag manifold of oriented
lines within oriented 4-planes in $\R^{8}$, similarly to (\cite{GZ}).
This gives rise to the following diagram:\[
\begin{array}{ccc}
G_{o}(3,S^{7}) & \longrightarrow & F_{o}(1,4,\R^{8})\\
h\Big\uparrow\pi\Big\downarrow &  & \pi\Big\downarrow\\
S^{7} &  & G(4,\R^{8})\end{array}.\]
$G_{o}(4,8)$ has two universal 4-plane bundles, $U(4,8)$ and $V(4,8)$,
defined by \begin{eqnarray*}
U(4,8) & := & \cup_{x\in G_{o}(4,8)}x\\
V(4,8) & := & \cup_{x\in G_{o}(4,8)}x^{\perp.}\end{eqnarray*}
The respective Euler classes $E(U)$ and $E(V)$ satisfy $E(U)\cup E(V)=0$
in $H^{8}(G_{o}(4,8))$. Similarly, the respective first Pontryagin
classes $P_{1}(U)$and $P_{1}(V)$ satisfy the same relationship,
$P_{1}(U)\cup P_{1}(V)=0$. In particular, if $\omega$ is the universal
connection on $U(4,8)$ defined by Narasimhan and Ramanan (cf. \cite{NR}),
and $\omega^{*}$ is the {}``dual'' connection on $V(4,8)$, then
the associated Euler forms, $e(\Omega)$ and $e(\Omega^{*})$, satisfy
$e(\Omega)\wedge e(\Omega^{*})=0$ (respectively, the first Pontryagin
forms). Then, consider the form \[
\Phi:=C\, Te(\omega)\wedge e(\Omega^{*}),\]
which is well-defined on $F_{o}(1,4,\R^{8})$ as well as on the frame
bundle $FU(4,8)$ of oriented orthonormal frames on $U(4,8)$, which
is an $SO(3)$-principal bundle over $F_{o}(1,4,\R^{8})$. Here, $Te(\omega)$
is the transgressive Chern-Simons form corresponding to the Euler
form $e(\Omega)$ of $U(4,8)$ \cite{CS}. Because $d(Te(\omega))=e(\Omega)$
(again, either as a form on the frame bundle, or on the associated
bundle $F_{o}(1,4,\R^{8})$), we have that $d\Phi=0$. The constant
$C$ is simply chosen so that the comass of $\Phi$ is one.

That $\Phi$ is well-defined on $F_{o}(1,4,\R^{8})$ is perhaps not
obvious. However, the original version of the transgressive form $Te(\omega)$
was defined by Chern on the sphere bundle, not the frame bundle \cite{C}.
That same construction applies here. When restricted to vertical directions,
those tangent to the 3-sphere fiber of $F_{o}(1,4,\R^{8})\rightarrow G_{o}(4,8)$,
$Te(\omega)$ is the volume form of the fibers.

\subsection{Calculations}

We will consider $G_{o}(4,8)$ as $SO(8)/SO(4)\times SO(4)$, and
the principal bundle $FU(4,8)$ as $SO(8)/I\times SO(4)$. The universal
connection $\omega$ on $FU(4,8)$ can be defined as the truncation
of the restriction of the Maurer-Cartan form on $o(8)$, denoted $\mu=\left[\mu_{ij}\right]$,
to the tangents to $FU(4,8)$. That is, the components of the connection
$\omega_{ij}$ are defined for $i,j\in\{1,\ldots,4\}$ and $\omega_{ij}(A)=A_{ij}$
for any \[
A\in T_{*}(U(4,8),(U_{0},\left\{ e_{1},\ldots,e_{4}\right\} ))=\left\{ \left.A\in o(8)\right|A=\left[\begin{array}{cc}
R & S\\
-S^{t} & 0\end{array}\right],\, R\in o(4)\right\} ,\]
if $U_{0}=\R^{4}\times0\subset\R^{8}$, with basis $\left\{ e_{1},\ldots,e_{4}\right\} $.
By homogeneity, all calculations in $FU(4,8)$ can be taken to be
at this point. 

The curvature $\Omega$ of this connection is given by $\Omega_{ij}(X,Y)=-\omega_{ij}([X,Y])$
for left-invariant vector fields that are horizontal at $U_{0}$,
that is, of the form $\left[\begin{array}{cc}
0 & S\\
-S^{t} & 0\end{array}\right]$. In terms of the Maurer-Cartan form, $\Omega_{ij}=+\sum_{k=5}^{8}\mu_{ik}\wedge\mu_{jk}$,
for $i,j\in\{1,\ldots,4\}$. 

Similarly, the connection $\omega^{*}$ on the dual principal bundle
$FV(4,8)=SO(8)/SO(4)\times I$ at $U_{0}=\R^{4}\times0\subset\R^{8}$
is the restriction of the same Maurer-Cartan form $\mu$ to the other
$4\times4$ block, and the curvature $\Omega_{kl}^{*}=\sum_{i=1}^{4}\mu_{ik}\wedge\mu_{il}$,
for $k,l\in\left\{ 5,\ldots,8\right\} $. Either of the tangent spaces
to these principal bundles can be canonically embedded into the tangent
space $o(8)$ of $SO(8)$ at the identity.

The Euler form $e(\Omega)$ of $FU(4,8)$ is the form \begin{eqnarray*}
e(\Omega) & := & \frac{1}{2\pi^{2}}\left(\Omega_{12}\wedge\Omega_{34}-\Omega_{13}\wedge\Omega_{24}+\Omega_{14}\wedge\Omega_{23}\right)\\
 & = & \frac{1}{2\pi^{2}}\left(\mu_{1k}\wedge\mu_{2k}\wedge\mu_{3l}\wedge\mu_{4l}-\mu_{1k}\wedge\mu_{3k}\wedge\mu_{2l}\wedge\mu_{4l}+\mu_{1k}\wedge\mu_{4k}\wedge\mu_{2l}\wedge\mu_{3l}\right),\end{eqnarray*}
where the sum is taken over all $k,l\in\left\{ 5,\ldots,8\right\} $.
Dually, the Euler form $e(\Omega^{*})$ of $FV(4,8)$ is the form
\begin{eqnarray*}
e(\Omega^{*}) & := & \frac{1}{2\pi^{2}}\left(\Omega_{56}\wedge\Omega_{78}-\Omega_{57}\wedge\Omega_{68}+\Omega_{58}\wedge\Omega_{67}\right)\\
 & = & \frac{1}{2\pi^{2}}\left(\mu_{5i}\wedge\mu_{6i}\wedge\mu_{7j}\wedge\mu_{8j}-\mu_{5i}\wedge\mu_{7i}\wedge\mu_{6j}\wedge\mu_{8j}+\mu_{5i}\wedge\mu_{8i}\wedge\mu_{6j}\wedge\mu_{7j}\right),\end{eqnarray*}
where the sum is taken over all $i,j\in\left\{ 1,\ldots,4\right\} $.

\begin{prop}
$e(\Omega)\wedge e(\Omega^{*})\equiv0$.
\end{prop}
\begin{proof}
Each monomial in this product is of the form \[
\mu_{1k}\wedge\mu_{2k}\wedge\mu_{3l}\wedge\mu_{4l}\wedge\mu_{5i}\wedge\mu_{6i}\wedge\mu_{7j}\wedge\mu_{8j}\]
 or a permutation thereof. $k$ can be either $5,\,6,\,7$ or $8$.
If $k$ is, say, $5,$ then $i$ cannot be $1$ or $2,$ thus must
be $i=3\,\textrm{or}\,4$. Thus $l\neq5,6$, so $l=7\,\textrm{or}\,8$,
and finally, $j=1\,\textrm{or}\,2$. No matter which choices are made,
two of the indices between $1$ and $4$ will occur once, and the
other two will occur three times, and similarly for the indices from
$5$ to $8$. Thus, each monomial is determined by the multi-indices
that occur with one index singly. For example, $2,\,4,\,6,\,\textrm{and}\,8$
occur singly, paired as $25,\,47,\,36,\,\textrm{and}\,18$ in exactly
two terms, \begin{eqnarray*}
 &  & +\mu_{15}\wedge\mu_{25}\wedge\mu_{37}\wedge\mu_{47}\wedge\mu_{53}\wedge\mu_{63}\wedge\mu_{71}\wedge\mu_{81},\,\textrm{and}\\
 &  & +\mu_{17}\wedge\mu_{47}\wedge\mu_{25}\wedge\mu_{35}\wedge\mu_{51}\wedge\mu_{81}\wedge\mu_{63}\wedge\mu_{73}.\end{eqnarray*}
However, using the fact that $\mu_{ik}=-\mu_{ki}$ and the exterior
product, these terms cancel. Since all terms are permutations of these,
all terms cancel in pairs.
\end{proof}
Thus, the form $\Phi:=C\, Te(\omega)\wedge e(\Omega^{*})$ is indeed
closed. That it is well-defined can be traced back to early versions
of the Chern-Simons theory, such as \cite{C}. Alternately, it can
be directly verified from the local expression for $Te(\omega)$ in
terms of the Maurer-Cartan form $\mu$. That is, as a form on $F_{o}(1,4,\R^{8})$,
at the point $x_{0}:=(e_{1},W)$, $e_{1}\in W=\R^{4}\times\{0\}\subset\R^{8}$,
since all the $\omega_{ij}$ tangent to $F_{o}(1,4,\R^{8})$ have
one of $i=1$ or $j=1$,\begin{eqnarray*}
Te(\omega) & := & \frac{1}{2\pi^{2}}\left(\omega_{12}\Omega_{34}-\omega_{13}\Omega_{24}+\omega_{14}\Omega_{23}\right.\\
 &  & \left.-\frac{1}{6}\left(\omega_{12}\left(\left[\omega,\omega\right]\right)_{34}-\omega_{13}\left(\left[\omega,\omega\right]\right)_{24}+\omega_{14}\left(\left[\omega,\omega\right]\right)_{23}\right)\right)\\
 & = & \frac{1}{2\pi^{2}}\left(\mu_{12}\wedge\mu_{3k}\wedge\mu_{4k}-\mu_{13}\wedge\mu_{2k}\wedge\mu_{4k}+\mu_{14}\wedge\mu_{2k}\wedge\mu_{3k}\right.\\
 &  & \left.+\frac{1}{3}\left(\mu_{12}\wedge\mu_{13}\wedge\mu_{14}-\mu_{13}\wedge\mu_{12}\wedge\mu_{14}+\mu_{14}\wedge\mu_{12}\wedge\mu_{13}\right)\right)\\
 & = & \frac{1}{2\pi^{2}}\left(\mu_{12}\wedge\mu_{13}\wedge\mu_{14}+\mu_{12}\wedge\mu_{3k}\wedge\mu_{4k}-\mu_{13}\wedge\mu_{2k}\wedge\mu_{4k}+\mu_{14}\wedge\mu_{2k}\wedge\mu_{3k}\right),\end{eqnarray*}
where the sum is over $k$ from $5$ to $8$. As a left-invariant
form on $SO(8)$, it is straightforward to see that it is invariant
under the adjoint action of the isotropy subgroup $1\times SO(3)\times I_{4}\subset SO(8)$,
so descends to a form on $F(1,4,\R^{8})$. 

It remains to find the maximum of $\Phi(W)$ for $7$-planes $W$
in the total space of $\pi\: F(1,4,\R^{8})\rightarrow G(4,8)$.

Not all vertical directions (those in the 3-sphere fiber of $\pi$)
and combinations within $\Lambda_{3}(T_{*}(F(1,4,\R^{8})))$ are detected
by the form $\Phi$ above. That is, let $W\subset T_{*}(F(1,4,\R^{8}))$
be a $7$-dimensional subspace. Then, for some basis of $T_{*}(F(1,4,\R^{8}))$,
with vertical directions $v_{1},\, v_{2},\, v_{3}$ and horizontal
basis $\left\{ e_{1},\ldots,e_{16}\right\} $, $W=(a_{1}v_{1}+h_{1})\wedge(a_{2}v_{2}+h_{2})\wedge(a_{3}v_{3}+h_{3})\wedge(e_{1}\wedge\cdots\wedge e_{4})$
as a unit element of $\Lambda_{7}\left(T_{*}(F(1,4,\R^{8}))\right)$.
$\left\Vert h_{i}\right\Vert =b_{i}=\sqrt{1-a_{i}^{2}}$. This simply
states that no more than 3 directions can have independent vertical
components. At the point $x_{0}$, if we denote $a_{i}=\cos(\theta_{i})$
and $b_{i}=\sin(\theta_{i})$,\begin{eqnarray*}
\Phi(W) & \leq & C\left|a_{1}a_{2}a_{3}+a_{1}b_{2}b_{3}-b_{1}a_{2}b_{3}+b_{1}b_{2}a_{3}\right|\left|e(\Omega^{*})(e_{1}\wedge\cdots\wedge e_{4})\right|\\
 & \leq & C\left|\cos(\theta_{1})\left(\cos(\theta_{2})\cos(\theta_{3})+\sin(\theta_{2})\sin(\theta_{3})\right)-\sin(\theta_{1})\left(\cos(\theta_{2})\sin(\theta_{3})-\sin(\theta_{2})\cos(\theta_{3})\right)\right|\cdot\\
 &  & \cdot\left|e(\Omega^{*})(e_{1}\wedge\cdots\wedge e_{4})\right|\\
 & = & C\left|\cos(\theta_{1}+(\theta_{2}-\theta_{3}))\right|\left|e(\Omega^{*})(e_{1}\wedge\cdots\wedge e_{4})\right|\\
 & \leq & C\left|e(\Omega^{*})(e_{1}\wedge\cdots\wedge e_{4})\right|,\end{eqnarray*}
Thus the maximum is achieved when (among other values) all three $\theta_{i}$
are 0, as long as the remaining vectors form a $4$-plane maximizing
$e(\Omega^{*})$. It is not clear whether other values of $\theta_{i}$
will achieve this maximum, since the mixed parts of $Te(\Omega)$
are only bounded by those values. However, the maximum is clearly
achieved when all $\theta_{i}=0$. 

Let $\pi_{U}:F_{o}(1,4,\R^{8})\rightarrow G_{o}(4,8)$ be the fibration
associated with the unit sphere bundle of the universal bundle $U$,
so that $(x,W)$, where $x\in W$ is a unit vector in the $4$-plane
$W$, is mapped to $\pi_{u}(x,W):=W\in G_{o}(4,8)$. The other fibration
$\pi_{V}:F_{o}(1,4,\R^{8})\rightarrow G_{o}(4,8)$, associated with
the dual bundle $V$, maps the same $(x,W)$ onto $\pi_{V}(x,W):=W^{\perp}$.
$e(\Omega^{*})$, as a form on $F_{o}(1,4,\R^{8})$, is the $\pi_{V}$-horizontal
lift of the form $e(\Omega^{*})$ on $G_{o}(4,8)$. That is clearly
maximized on some collection of $\pi_{V}$-horizontal 4-planes tangent
to $G_{o}(4,8)$ at $W$. $Te(\omega)$ is maximized on the 3-sphere
fibers of $W$, that is $\left\{ \left.(x,W)\right|x\in W\right\} $,
as described above. Since these two spaces are orthogonal, then \[
\Phi:=C\, Te(\omega)\wedge e(\Omega^{*})\]
 will be maximized on any $7$-plane which is the sum of a $\pi_{V}$-horizontal
lift of a $4$-plane maximizing $e(\Omega^{*})$ (perpendicular to
$W$) and the $3$-plane tangent to the unit sphere in $W$ at $x$.
However, not all $4$-planes orthogonal to $W$ will maximize $\Phi$. 

Since \begin{eqnarray*}
e(\Omega^{*}) & := & \frac{1}{2\pi^{2}}\left(\Omega_{56}\wedge\Omega_{78}-\Omega_{57}\wedge\Omega_{68}+\Omega_{58}\wedge\Omega_{67}\right)\\
 & = & \frac{1}{2\pi^{2}}\left(\mu_{5i}\wedge\mu_{6i}\wedge\mu_{7j}\wedge\mu_{8j}-\mu_{5i}\wedge\mu_{7i}\wedge\mu_{6j}\wedge\mu_{8j}+\mu_{5i}\wedge\mu_{8i}\wedge\mu_{6j}\wedge\mu_{7j}\right),\end{eqnarray*}
if $E_{ij}$ is the basis of tangent vectors dual to $\mu_{ij}$,
$e(\Omega^{*})(E_{15},E_{16},E_{17},E_{18})=3/2\pi^{2}$. It is straightforward
to see that $e(\Omega^{*})(E_{i5},E_{i6},E_{i7},E_{i8})=3/2\pi^{2}$
for any $i=1\ldots4$, and $e(\Omega^{*})(E_{i5},E_{i6},E_{j7},E_{j8})=1/2\pi^{2}$
for $i\neq j$, or, more generally, if $k_{1},\ldots,k_{4}$ are a
permutation of $5,\ldots,9$, then when $i\neq j$, $e(\Omega^{*})(E_{ik_{1}},E_{ik_{2}},E_{jk_{3}},E_{jk_{4}})=\pm1/2\pi^{2}$,
where the sign is the sign of the permutation. Finally, if $\{ i_{1},i_{2},i_{3},i_{4}\}$
consist of more than two distinct values (and not two pairs of values),
or if $\{ k_{1},k_{2},k_{3},k_{4}\}$ does not consist of some permutation
of $\{5,6,7,8\},$ then \[
e(\Omega^{*})(E_{i_{1}k_{1}},E_{i_{2}k_{2}},E_{i_{3}k_{3}},E_{i_{4}k_{4}})=0.\]

\begin{thm}
The singular foliation $NS$ of $S^{7}$ consisting of all great 3-spheres
containing a common great 2-sphere minimizes volume of all three-dimensional
singular foliations on $S^{7}$ with that singular locus and limiting
behavior, and provides a lower bound for the volume of all regular
three-dimensional oriented foliations of $S^{7}$.
\end{thm}
\begin{proof}
For any decomposable, unit $\xi\in\Lambda_{4}(G(4,8),W_{0})$ which
is tangent to the variety $G(4,8)$ at $W_{0}$, \[
\xi=\sum_{i_{1},\cdots,i_{4},k_{1}\leq\cdots\leq k_{4}}\xi_{i_{1},\ldots,i_{4},k_{1},\ldots,k_{4}}E_{i_{1}k_{1}}\wedge\cdots\wedge E_{i_{4}k_{4}}.\]
Since $\xi$ is decomposable, $\xi$ satisfies the Pl\"{u}cker condition
$\xi\wedge\xi=0$, implying that, in particular (restricting to the
case where $\{ k_{1},\ldots,k_{4}\}=\{5,6,7,8\}$ since otherwise
$e(\Omega^{*})=0$), and denoting $\xi_{i,j,k,l,5,6,7,8}$ by $\xi_{i,j,k,l}$,
\begin{eqnarray*}
\xi_{1,1,1,1}\xi_{2,2,2,2}-\xi_{2,1,1,1}\xi_{1,2,2,2}-\xi_{1,2,1,1}\xi_{2,1,2,2}-\xi_{1,1,2,1}\xi_{2,2,1,2}\\
-\xi_{1,1,1,2}\xi_{2,2,2,1}+\xi_{1,1,2,2}\xi_{2,2,1,1}+\xi_{1,2,1,2}\xi_{2,1,2,1}+\xi_{1,2,2,1}\xi_{2,1,1,2} & = & 0.\end{eqnarray*}
and similarly for all other such combinations. Thus, \begin{eqnarray*}
 &  & (\xi_{1,1,1,1}+\xi_{2,2,2,2}+\xi_{3,3,3,3}+\xi_{4,4,4,4})^{2}\\
 & = & \xi_{1,1,1,1}^{2}+\xi_{2,2,2,2}^{2}+2\xi_{1,1,1,1}\xi_{2,2,2,2}+\cdots\\
 & = & \xi_{1,1,1,1}^{2}+\xi_{2,2,2,2}^{2}+2\xi_{2,1,1,1}\xi_{1,2,2,2}+2\xi_{1,2,1,1}\xi_{2,1,2,2}+2\xi_{1,1,2,1}\xi_{2,2,1,2}\\
 &  & +2\xi_{1,1,1,2}\xi_{2,2,2,1}-2\xi_{1,1,2,2}\xi_{2,2,1,1}-2\xi_{1,2,1,2}\xi_{2,1,2,1}-2\xi_{1,2,2,1}\xi_{2,1,1,2}+\cdots\\
 & \leq & \xi_{1,1,1,1}^{2}+\xi_{2,2,2,2}^{2}+\xi_{2,1,1,1}^{2}+\xi_{1,2,2,2}^{2}+\xi_{1,2,1,1}^{2}+\xi_{2,1,2,2}^{2}+\xi_{1,1,2,1}^{2}+\xi_{2,2,1,2}^{2}\\
 &  & +\xi_{1,1,1,2}^{2}+\xi_{2,2,2,1}^{2}+\xi_{1,1,2,2}^{2}+\xi_{2,2,1,1}^{2}+\xi_{1,2,1,2}^{2}+\xi_{2,1,2,1}^{2}\\
 &  & +\xi_{1,2,2,1}^{2}+\xi_{2,1,1,2}^{2}+\cdots\\
 & \leq & 1,\end{eqnarray*}
since $\xi$ is a unit. Thus, on any such $\xi$, \[
e(\Omega^{*})(\xi)\leq3/2\pi^{2},\]
the maximum being achieved on those $\xi$ so that $(\xi_{1,1,1,1,5,6,7,8}+\xi_{2,2,2,2,5,6,7,8}+\xi_{3,3,3,3,5,6,7,8}+\xi_{4,4,4,4,5,6,7,8})=1.$
Those 4-planes, in contrast to the complex case studied by Gluck and
Ziller, are \emph{not} those which are tangent 4-planes in $T_{*}(G(4,8),W_{0})$
to the quaternionic projective space ${\mathbb{HP}}^{1}$ under any
quaternionic structure on $\R^{8}$ for which $W_{0}$ is quaternionic.
Those 4-planes can be easily shown to evaluate to half the maximum
possible value. 

In fact, if $M$ is the manifold \[
M:=\left\{ \left.x\wedge e_{2}\wedge e_{3}\wedge e_{4}\right|x\perp\{ e_{2},e_{3},e_{4}\},\,\left\Vert x\right\Vert =1\right\} ,\]
then $M\cong S^{4}$, and the corresponding {}``foliation'' on $S^{7}$
will be singular. The tangent planes to $M$ at each point clearly
maximize the value of $e(\Omega^{*})$. The corresponding singular
foliation on $S^{7}$ is the set of all great 3-spheres that are intersections
of $S^{7}$ with a plane $W=span\{ x,e_{2},e_{3},e_{4}\}$ for some
unit $x\perp\{ e_{2},e_{3},e_{4}\}$, which is singular on the $S^{2}$
common to all leaves. However, this singular distribution will indeed
be calibrated by this form, so minimizes volume among, at least, all
singular foliations with the same singular set; in this case, a totally-geodesic
$S^{2}$ which is the intersection of any two leaves of the foliation. 

As with the case for one-dimensional leaves, this singular foliation
actually corresponds to half of the current $\pi^{-1}(M)\subset F_{0}(1,4,8)\cong G(3,S^{7})$,
since the leaf corresponding to the 4-plane $x\wedge e_{2}\wedge e_{3}\wedge e_{4}$
is the same set as that leaf corresponding to $(-x)\wedge e_{2}\wedge e_{3}\wedge e_{4}$
with the opposite orientation. The 3-plane common to all 4-planes
separates each into two half-spaces. Choose the half-space consistent
with a chosen orientation on the common 3-plane, which then restricts
the fibers to hemispheres which still provides a singular foliation
of $S^{7}$. Since this (non-cycle) current $NS\subset G_{o}(3,S^{7})$
has boundary $S^{2}\times S^{4}\subset\left.G_{o}(3,S^{7})\right|_{S^{2}}\cong S^{2}\times G_{o}(3,7),$
which is not itself a boundary, $NS$ does not extend to a cycle.
Thus, the fact that $\Phi$ calibrated $NS$ only implies that $NS$
represents a singular foliation on $S^{7}$ which is volume minimizing
among foliations with the same singular locus.

However, similarly to \cite{GZ}, it follows that the full preimage
$S:=\pi^{-1}(M)$, which is also calibrated by $\Phi$ and is a cycle,
minimizes mass among currents homologous to twice the homology class
of a foliation (all foliations by 3-manifolds are homologous as maps
into $G_{o}(3,S^{7})$). If there were a (singular or regular) volume-minimizing
foliation represented by a cycle $C$, then the mass of $2C$ could
not be less than the mass of $S$, so that the mass of $NS$ does
represent a lower bound of volumes of foliations of dimension $3$
on $S^{7}$. 

\end{proof}
It remains an open question whether the Hopf fibration minimizes volume
among 3-dimensional regular foliations of $S^{7}$. However, the Hopf
fibration (a regular foliation) does have twice the volume of the
singular foliation $NS$.

\section{Generalizations}

It is a straightforward generalization of these computations to show
that the corresponding sphere $M$ maximizes the corresponding form
$e(\Omega^{*})$ in $G_{o}(4,4n+4)$, showing that similar singular
foliations by 3-manifolds minimize volume among all (singular) foliations
of $S^{4n+3}$ with the given singular set. 

\begin{thm}
The singular foliation of $S^{4n+3}$ consisting of all great 3-spheres
containing a common great 2-sphere minimizes volume of all three-dimensional
singular foliations on $S^{4n+3}$ with that singular locus and limiting
behavior, and provides a lower bound for the volume of all regular
three-dimensional oriented foliations of $S^{4n+3}$.
\end{thm}
Similarly, the same methods will show that the Hopf fibration of $S^{15}$
by great $7$-spheres, the fibers of the Cayley projective plane,
the fibers of the fibration \[
\begin{array}{ccc}
S^{7} & \rightarrow & S^{15}\\
 &  & \downarrow\\
 &  & S^{8}\end{array},\]
 will not minimize volume among all singular foliations of that space
as well, but rather the {}``longitudes'', great 7-spheres foliating
$S^{15}$ except for a great 6-sphere common to all leaves, will be
a volume-minimizing singular foliation. 

\begin{thm}
The singular foliation $NS $ of $S^{15}$ consisting of all great
7-spheres containing a common great 6-sphere minimizes volume of all
7-dimensional singular foliations on $S^{15}$ with that singular
locus and limiting behavior, and provides a lower bound for the volume
of all regular three-dimensional oriented foliations of $S^{15}$.
\end{thm}

\end{document}